\documentclass[10pt,twoside]{article}
\usepackage{Latex-document}
\almvol{00}
\almttone{Frontiers of Science Awards}
\almtttwo{for Math/TCIS/Phys}
\firstpage{1}
\usepackage[english]{babel}
\usepackage[utf8]{inputenc}

\usepackage{multicol}

\usepackage{tabu}

\usepackage{amsmath}
\usepackage{amssymb}
\usepackage{graphics}

\numberwithin{equation}{section}

\newcommand{\LL}{{\mathbb{L}}}

\newcommand{\Q}{\mathbb{Q}}
\newcommand{\Z}{\mathbb{Z}}

\newcommand{\rarr}{\rightarrow}

\newcommand{\Mb}{\overline{\mathcal{M}}}

\newcommand{\com}{{\mathbb C}}

\newcommand{\bZ}{\mathsf{Z}}

\newcommand{\PC}{{\mathcal{P}}}

\newcommand\FF{\mathbb F}
\newcommand\ZZ{\mathsf Z}
\newcommand{\oh}{{\mathcal O}}

\newcommand{\OO}{{\mathcal{O}}}

\newcommand{\Hilb}{\mathsf{Hilb}^n(\mathbb{C}^2)}
\newcommand{\Sym}{\mathsf{Sym}^n}

\begin{document}

\markboth{\hfill{\rm Rahul Pandharipande} \hfill}{\hfill
  {\rm Moduli of curves and moduli of sheaves \hfill}}

\title{Moduli of curves and moduli of sheaves}

\author{Rahul Pandharipande}

\begin{abstract}

  Relationships between moduli spaces of curves and sheaves
  on $3$-folds are presented starting with the
  Gromov-Witten/Donaldson-Thomas correspondence proposed more than 20 years ago with D. Maulik, N. Nekrasov, and A. Okounkov.
  The descendent and relative correspondences as developed with A. Pixton
  in the context of stable pairs led to the proof of the 
  correspondence for the Calabi-Yau quintic 3-fold. More recently, the study of correspondences in families has played an
  important role in connection with other basic moduli problems in algebraic geometry. The full conjectural framework is presented here in the context
  of families of $3$-folds.

\end{abstract}

\maketitle

\setcounter{tocdepth}{1}
\tableofcontents

\section{Introduction}
Let $X$ be a nonsingular complex projective $3$-fold. The counting
of algebraic curves in $X$ can be approached by stable maps in
Gromov-Witten theory and by stable pairs in Donaldson-Thomas theory.
In a series of conjectures starting with the original
Gromov-Witten/Donaldson-Thomas correspondence \cite{MNOP1,MNOP2}
for ideal sheaves
formulated with D. Maulik, N. Nekrasov, and A. Okounkov,
an equivalence is proposed between the 
different approaches to counting via non-trivial
transformations.
My goal here is to present the Gromov-Witten/Pairs
descendent correspondence
in a general form for families and to discuss a few cases in which the
correspondence has been proven (including the case of the Calabi-Yau
quintic $3$-fold established with A. Pixton in \cite{PPQ}).

The full Gromov-Witten/Pairs descendent correspondence
for families in the form of Conjectures II and III of Section 5
is more general than in the previous
perspective on the subject presented in Conjecture I of Section 3. Sections 2-4 of paper are essentially
a survey of the conjectures and results of \cite{PPDC,PPQ} and related work.
Other introductions to the subject can be found in \cite{PanDesc,13H}.

\section{Stable maps and stable pairs}

\subsection{Families}
We will state all the conjectures in the context of families.
Let $$\nu:\mathcal{X} \rightarrow \mathcal{Y}$$ be a flat projective
family of
nonsingular complex $3$-folds over an irreducible quasi-projective
base scheme $\mathcal{Y}$. Later, in the context of relative and
log theories, the morphism $\nu$ may be taken to be log smooth instead
of smooth. The focus of the paper, however, will be on the smooth case.

\subsection{Gromov-Witten theory: stable maps}

Let
 $\Mb_{g,r}(X,\beta)$ denote the moduli space of
 $r$-pointed stable maps from connected genus $g$ curves to
a nonsingular projective $3$-fold
 $X$ representing a
 class  $\beta\in H_2(X,\Z)$. Let 
$$\text{ev}_i: \Mb_{g,r}(X,\beta) \rarr X\, , \ \ \  \LL_i \rarr \Mb_{g,r}(X,\beta)$$
denote the evaluation maps and the cotangent line bundles associated to
the marked points.

We will be interested in the moduli space of stable
maps 
$$\Mb_{g,r}(\nu,\beta) \ \ \ \text{associated to} \ \ \ 
\nu:\mathcal{X} \rightarrow \mathcal{Y}\, ,$$
where $\beta$ is a fiber class well-defined for the family.
The points of $\Mb_{g,r}(\nu,\beta)$ are stable maps to
the fibers of $\nu$.

Let 
$\gamma_1, \ldots, \gamma_r\in H^*(\mathcal{X})$
be cohomology{\footnote{All homology
and cohomology groups will be taken with $\Q$-coefficients
unless explicitly denoted otherwise.}} classes, and
let $$\psi_i = c_1(\LL_i) \in H^2(\Mb_{g,n}(\nu,\beta))\, .$$
The {\em standard descendent} insertions, denoted by $\tau_k(\gamma)$, correspond 
to the classes $\psi_i^k \text{ev}_i^*(\gamma)$ on the moduli space
of stable maps. 
Let
$$\Big\langle \tau_{k_1}(\gamma_{1}) \cdots
\tau_{k_r}(\gamma_{r})\Big\rangle_{g,\beta}^{\mathsf{GW}} =
\epsilon_{\mathsf{GW}*} \left(\prod_{i=1}^r \psi_i^{k_i} \text{ev}_i^*(\gamma_{i})
\,  \cap \, [\Mb_{g,r}(\nu,\beta)]^{vir}
\right) \in H_*(\mathcal{Y})$$
denote the descendent
Gromov-Witten invariants where
$$\epsilon_{\mathsf{GW}}: \Mb_{g,r}(\nu,\beta) \rightarrow \mathcal{Y}\, $$
is the map to the base of the family.
Foundational aspects of deformation theory
and the virtual class
are treated in \cite{BehFan, LiTian}.

Let $C$ be a possibly disconnected curve with at worst nodal singularities.
The genus of $C$ is defined by $1-\chi(\oh_C)$. 
Let $\Mb'_{g,r}(\nu,\beta)$ denote
the moduli space of maps mapping to fibers of $\nu$
with possibly {disconnected} domain
curves $C$ of genus $g$ with {\em no} collapsed connected components.
In particular, 
$C$ must represent a {nonzero} fiber class $\beta$.

We define the descendent invariants in the disconnected 
case by
$$\Big\langle \tau_{k_1}(\gamma_{1}) \cdots
\tau_{k_r}(\gamma_{r})\Big\rangle^{\mathsf{GW}'}_{g,\beta} =
\epsilon'_{\mathsf{GW}*} \left(\prod_{i=1}^r \psi_i^{k_i} \text{ev}_i^*(\gamma_{i})
\,  \cap \, [\Mb'_{g,r}(\nu,\beta)]^{vir}
\right) \in H_*(\mathcal{Y})\, ,$$
where $\epsilon'_{\mathsf{GW}}: \Mb'_{g,r}(\nu,\beta) \rightarrow \mathcal{Y}$.
The associated Gromov-Witten descendent partition function is defined by
\begin{equation}
\label{abcmap}
\bZ'_{\mathsf{GW}}\Big(\nu;u\ \Big|\ \prod_{i=1}^r \tau_{k_i}(\gamma_{i})\Big)_\beta = 
\sum_{g\in{\mathbb Z}} \Big \langle \prod_{i=1}^r
\tau_{k_i}(\gamma_{i}) \Big \rangle^{\mathsf{GW}'}_{g,\beta} \ u^{2g-2}.
\end{equation}
Since the domain components must map non-trivially, an elementary
argument shows the genus $g$ in the  sum on the right side of \eqref{abcmap} is bounded from below. 
As a result, the partition function \eqref{abcmap} is a Laurent series
in $u$ with coefficients in $H_*(\mathcal{Y})$.

\subsection{Product evaluations for stable maps}

\label{prgw}

Let $\nu^r: \mathcal{X}^r \rightarrow \mathcal{Y}$ be the fiber
product over $\mathcal{Y}$,
$$ \mathcal{X}^r = \underbrace{\mathcal{X} \times_\mathcal{Y}  \cdots
  \times_{\mathcal{Y}} \mathcal{X}}_{r}\,.$$
We can consider cohomology classes
$\delta\in H^*(\mathcal{X}^r)$ and
define {\em general descendent} insertions
$$
\tau_{k_1,\ldots,k_r}(\delta)= \prod_{i=1}^r \psi_i^{k_i} \cdot \text{ev}_{\text{full}}^*(\delta)$$
using the evaluation map
$$\text{ev}_{\text{full}}\colon \Mb'_{g,r}(\nu,\beta)\to \mathcal{X}^r$$
obtained from the full  marking set.

We define the associated descendent invariants in the disconnected 
case by
$$\Big\langle \tau_{k_1,\ldots, k_r}(\delta) \Big\rangle^{\mathsf{GW}'}_{g,\beta} =
\epsilon'_{\mathsf{GW}*} \left(\prod_{i=1}^r \psi_i^{k_i} \cdot
  \text{ev}_{\text{full}}^*(\delta)
\,  \cap \, [\Mb'_{g,r}(\nu,\beta)]^{vir}
\right) \in H_*(\mathcal{Y})\, .$$
The associated Gromov-Witten descendent partition function is defined by
\begin{equation*}
\bZ'_{\mathsf{GW}}\Big(\nu;u\ \Big|\
                     \tau_{k_1,\ldots, k_r}(\delta)
\Big)_\beta = 
\sum_{g\in{\mathbb Z}} \Big \langle 
\tau_{k_1,\ldots, k_r}(\delta)
\Big \rangle^{\mathsf{GW}'}_{g,\beta} \ u^{2g-2}.
\end{equation*}

If  we have a factorization $\delta =\gamma_1 \otimes \cdots \otimes \gamma_r$ for classes
$\gamma_i \in H^*(\mathcal{X})$ pulled-back from the
$r$ distinct factor projections
$$\mathcal{X}^r \to \mathcal{X}\, ,$$
then 
$\tau_{k_1,\ldots,k_r}(\delta)= \prod_{i=1}^r \tau_{k_i}(\gamma_i)$,
and we have an equality of partition functions
$$\bZ'_{\mathsf{GW}}\Big(\nu;u\ \Big|\
                     \tau_{k_1,\ldots, k_r}(\delta)
                     \Big)_\beta =
                   \bZ'_{\mathsf{GW}}\Big(\nu;u\ \Big|\
                     \prod_{i=1}^r \tau_{k_i}(\gamma_i)
                     \Big)_\beta \, .$$
   However, $H^*(\mathcal{X}^r)$ is not in general spanned by
   such factorizations.{\footnote{If $\mathcal{Y}$ is
       a point, then the K\"unneth decomposition proves that factorizations
       span.}}

\subsection{Donaldson-Thomas theory: stable pairs}

\label{ddtt}

Let $X$ be a nonsingular complex projective $3$-fold.
A {\em{stable pair}} $(F,s)$ on $X$ is a coherent sheaf $F$ on $X$
 and a  section $s\in H^0(X,F)$ satisfying the following stability conditions:
\begin{itemize}
\item $F$ is \emph{pure} of dimension 1,
\item the section $s:\OO_X\to F$ has cokernel of dimensional 0.
\end{itemize}
To a stable pair, we associate the Euler characteristic and
the class of the support $C$ of $F$,
$$\chi(F)=n\in \mathbb{Z} \  \ \ \text{and} \ \ \ [C]=\beta\in H_2(X,\mathbb{Z})\,.$$
For fixed $n$ and $\beta$,
there is a projective moduli space of stable pairs $\PC_n(X,\beta)$. 

The moduli space of stable pairs is also defined in families of $3$-folds. 
We will be interested in the moduli space of stable
pairs
$$\PC_n(\nu,\beta) \ \ \ \text{associated to} \ \ \ 
\nu:\mathcal{X} \rightarrow \mathcal{Y}\, .$$
A foundational
treatment of the moduli space
of stable pairs 
is presented in \cite{pt} via the results of  Le Potier \cite{LePot}.
For a fixed fiber class $\beta$,
the moduli space $\PC_n(\nu,\beta)$ is
empty for all sufficiently negative $n$.

Denote the universal stable pair over $\mathcal{X}\times_{\mathcal{Y}} \PC_{n}(\nu,\beta)$ by
$$\OO_{\mathcal{X}\times_{\mathcal{Y}} \PC_n(\nu,\beta)} \stackrel{s\ }{\rightarrow}
\FF\, .$$
For $y\in \mathcal{Y}$ and
a stable pair $(F,s)\in \PC_{n}(\mathcal{X}_y,\beta)$, the restriction of
the universal stable pair
to the fiber
 $$\mathcal{X}_y \times (F,s) \ \subset\  
\mathcal{X}_y\times \PC_{n}(\mathcal{X}_y,\beta)
$$
is canonically isomorphic to $\OO_{\mathcal{X}_y}\stackrel{s\ }{\rightarrow} F$.
Let
$$\pi_{\mathcal{X}}\colon \mathcal{X}\times_{\mathcal{Y}} \PC_{n}(\nu,\beta)\to \mathcal{X},$$
$$\pi_\PC\colon \mathcal{X}\times_{\mathcal{Y}} \PC_{n}(\nu,\beta)
\to \PC_{n}(\nu,\beta)$$
 be the projections onto the first and second factors.
Since $\nu$ is smooth
and
$\FF$ is $\pi_\PC$-flat, $\FF$ has a finite resolution 
by locally free sheaves.
Hence, the Chern character of the universal sheaf $\FF$ on 
$\mathcal{X} \times_{\mathcal{Y}} \PC_n(\nu,\beta)$ is well-defined.

For each cohomology
class $\gamma\in H^*(\mathcal{X})$ and integer $k\in \mathbb{Z}_{\geq 0}$,
the action of the {\em standard descendent} insertion $\tau_k(\gamma)$ is defined by
$$
\tau_k(\gamma)=\pi_{\PC*}\big(\pi_{\mathcal{X}}^*(\gamma)\cdot \text{ch}_{2+k}(\FF)
\cap \pi_\PC^*(\ \cdot\ )\big)\, .$$ 

\noindent The pull-back
$\pi^*_\PC$ is well-defined in homology since $\pi_\PC$ is flat \cite{Ast}.

We may view the descendent action as defining a
cohomology class 
$$\tau_k(\gamma)\in H^*(\PC_n(\nu,\beta))$$
or as defining an endomorphism
$$\tau_k(\gamma):
H_*(\PC_{n}(\nu,\beta))\to H_*(\PC_{n}(\nu,\beta))\, .
$$

We define the stable pairs invariant with descendent insertions by
\begin{equation*}
\Big\langle \tau_{k_1}(\gamma_1)\ldots \tau_{k_r}(\gamma_r)
\Big\rangle_{\!n,\beta}^{\mathsf{P}}  = 
\epsilon_{\mathsf{P}*}
\left(\prod_{i=1}^r \tau_{k_i}(\gamma_i)\, \cap \, [\PC_n(\nu,\beta)]^{vir}\right)
\in H_*(\mathcal{Y})\, ,
\end{equation*}
where $\epsilon_{\mathsf{P}}: \PC_n(\nu,\beta) \rightarrow \mathcal{Y}$.
The partition function is 
\begin{equation*}
\ZZ_{\mathsf{P}}\Big(\nu;q\ \Big|   \prod_{i=1}^r \tau_{k_i}(\gamma_{i})
\Big)_\beta
=\sum_{n\in \mathbb{Z}} 
\Big\langle \prod_{i=1}^r \tau_{k_i}(\gamma_{i}) 
\Big\rangle_{\!n,\beta}^{\mathsf{P}} q^n.
\end{equation*}
Since $\PC_n(\nu,\beta)$ is empty for sufficiently negative
$n$, the partition function 
is a Laurent series in $q$ with coefficients in $H_*(\mathcal{Y})$.

\subsection{Product evaluations for stable pairs}

\label{prdt}

Let $\nu^r: \mathcal{X}^r \rightarrow \mathcal{Y}$ be the fiber
product over $\mathcal{Y}$ as in Section \ref{prgw}.
We can consider cohomology classes
$\delta\in H^*(\mathcal{X}^r)$ and define {\em general descendent}
  insertions for stable pairs by
$$
\tau_{k_1,\ldots,k_r}(\delta)=\pi_{\PC*}\big(\pi_{\mathcal{X}^r}^*(\delta)\cdot
\prod_{i=1}^r
\text{ch}_{2+k_i}(\FF_i)
\cap \pi_\PC^*(\ \cdot\ )\big)\,:
H_*(\PC_{n}(\nu,\beta))\to H_*(\PC_{n}(\nu,\beta))
$$
using the maps 
$$\pi_{\mathcal{X}^r}\colon \mathcal{X}^r\times_{\mathcal{Y}} \PC_{n}(\nu,\beta)\to \mathcal{X}^r,$$
$$\pi_\PC\colon \mathcal{X}^r\times_{\mathcal{Y}} \PC_{n}(\nu,\beta)
\to \PC_{n}(\nu,\beta)\,. $$
Here, $\mathbb{F}_i$ denotes the universal sheaf from the $i^{th}$ factor
$\mathcal{X}\times_{\mathcal{Y}} \PC_{n}(\mathcal{X},\beta)$.

We define the associated descendent invariants by
\begin{equation*}
\Big\langle \tau_{k_1,\ldots,k_r}(\delta)
\Big\rangle_{\!n,\beta}^{\mathsf{P}}  = 
\epsilon_{\mathsf{P}*}
\left(\tau_{k_1,\ldots,k_r}(\delta)\, \cap \, [\PC_n(\nu,\beta)]^{vir}\right)
\in H_*(\mathcal{Y})\, .
\end{equation*}
The partition function is 
\begin{equation*}
  \ZZ_{\mathsf{P}}\Big(\nu;q\ \Big|
\tau_{k_1,\ldots,k_r}(\delta)
  \Big)_\beta
=\sum_{n\in \mathbb{Z}} 
\Big\langle
\tau_{k_1,\ldots,k_r}(\delta)
\Big\rangle_{\!n,\beta}^{\mathsf{P}} q^n.
\end{equation*}
If we have a factorization $\delta =\gamma_1 \otimes \cdots \otimes \gamma_r$ for classes
$\gamma_i \in H^*(\mathcal{X})$ pulled-back from the
$r$ distinct factor projections, then
$$\ZZ_{\mathsf{P}}\Big(\nu;q\ \Big|   \tau_{k_1,\ldots,k_r}(\delta)
\Big)_\beta =
\ZZ_{\mathsf{P}}\Big(\nu;q\ \Big|   \prod_{i=1}^r \tau_{k_i}(\gamma_{i})
\Big)_\beta\, .$$

\subsection{Dimension constraints}

The stable maps and stable pairs
descendent series 
\begin{equation}\label{fhh6}
\bZ'_{\mathsf{GW}}\Big(\nu;u\ \Big|\ \prod_{i=1}^r \tau_{k_i}(\gamma_{i})\Big)_\beta\, , \ \ \ \
\bZ_{\mathsf{P}}\Big(\nu;q\ \Big|\ \prod_{i=1}^r \tau_{k_i}(\gamma_{i})\Big)_\beta
\end{equation}
for the family $\nu:\mathcal{X} \rightarrow \mathcal{Y}$
both satisfy  dimension constraints.

For $\gamma_i\in H^{e_i}(\mathcal{X})$, the (real) dimension of the descendent
theories are
$$\tau_{k_i}(\gamma_i)\in H^{{e_i}+2k_i}(\Mb'_{g,r}(\nu,\beta))\, , \ \ \ \
\tau_{k_i}(\gamma_i)\in H^{{e_i}+2k_i-2}(P_{n}(\nu,\beta))\, .$$
The virtual dimensions are
$$\text{dim}_{\mathbb{C}} \, [\Mb'_{g,r}(\nu,\beta)]^{vir} = \int_\beta c_1(T_\nu)\, +\, r\, +\, \dim_{\mathbb{C}} \mathcal{Y}\, , \ \ \ \
\text{dim}_{\mathbb{C}} \, [\PC_n(\nu,\beta)]^{vir} = \int_\beta c_1(T_\nu)
\, +\, \dim_{\mathbb{C}} \mathcal{Y}\, $$
where $T_\nu$ is the relative tangent bundle.
The coefficients of both series \eqref{fhh6} therefore take values
in the {\em same} homology group
$$H_{2\int_\beta c_1(T_\nu)+2\dim_{\mathbb{C}} \mathcal{Y} - \sum_{i=1}^r (e_i+2k_i-2)}(\mathcal{Y})\, .$$

Is there a relationship
between the Gromov-Witten and stable pairs descendent series \eqref{fhh6}?
The immediate issues are:
\begin{enumerate}
\item[(i)] the series
involve different moduli spaces and universal structures,
\item[(ii)] the variables $u$ and $q$ of
the two series are different.
\end{enumerate}
The {\em descendent correspondence}
proposes a precise relationship between the Gromov-Witten and stable pairs
descendent series, but only after a change of variables to
address (ii).

\section{Descendent correspondence}

\subsection{Descendent notation}
%

A partition  $\sigma=(\sigma_1,\ldots,
\sigma_{{\ell}(\sigma)})$
of positive size $|\sigma|$ 
and length $\ell(\sigma)$ consists of integers $\sigma_k$ satisfying
the properties
$$ \sigma_1\geq \ldots\geq 
\sigma_{\ell(\sigma)} > 0\, , \ \ \ \ 
\sum_{k=1}^{\ell(\sigma)}\sigma_k = |\sigma|\, .$$
For the family of $3$-folds
$\nu:\mathcal{X} \rightarrow \mathcal{Y}$, 
let 
$$\iota_\Delta:\Delta_{\ell(\sigma)}\rightarrow \mathcal{X}^{\ell(\sigma)}$$
denote the
inclusion of the fiberwise 
small diagonal{\footnote{The small diagonal 
is the set of points of each fiber $\mathcal{X}_y^{\ell(\sigma)}$ for which
the coordinates $(x_1,\ldots, x_{\hat{\ell}})$ are all equal
$x_i=x_j$.}}
in the product $\mathcal{X}^{\ell(\sigma)}$ over $\mathcal{Y}$.
For $\gamma\in H^*(\mathcal{X})$, 
we write $$\gamma\cdot \Delta_{\widehat{\ell}} =\iota_{\Delta*}(\gamma) \in
H^*(\mathcal{X}^{\ell(\sigma)})\, .$$
The {\em diagonal descendent} insertion{\footnote{We follow the notation of
    Sections \ref{prgw} for stable maps and \ref{prdt} for
    stable pairs.}}
\begin{equation}\label{ff33}
\tau_{\sigma}(\gamma) = \tau_{\sigma_1-1, \ldots,
  \sigma_{\ell(\sigma)}-1}(\gamma\cdot \Delta_{\ell(\sigma)})
\end{equation}
will play an important role in the descendent correspondence.

\vspace{10pt}

\noindent {\bf Warning.}
If $\sigma=(\sigma_1)$ has a single part, we recover the
standard descendent
$$\tau_\sigma(\gamma)= \tau_{\sigma_1-1}(\gamma)\, .$$
The shift by 1 is a confusing aspect of the definition (but is necessary
when using language of partitions to describe descendents).
When reading a formula with descendents, care must be taken to
determine whether the descendent subscript is a {\em partition} or an
{\em integer}.
If the subscript is a partition, then the descendent is a {\em diagonal descendent}.
If the subscript is an integer (or a list of integers),
then the descendent is either
a {\em standard} or a {\em general descendent}.



\subsection{Correspondence matrix} \label{corrmat}
 
A central result of \cite{PPDC} is
the construction of
a universal correspondence matrix $\widetilde{\mathsf{K}}$ 
indexed by partitions
$\alpha$ and $\widehat{\alpha}$ of positive size
with{\footnote{Here, $i^2=-1$.}}
$$\widetilde{\mathsf{K}}_{\alpha,\widehat{\alpha}}\in 
\mathbb{Q}[i,c_1,c_2,c_3]((u))\, . $$
The elements of $\widetilde{\mathsf{K}}$
are constructed from the capped descendent vertex \cite{PPDC}
and satisfy two basic properties:

\begin{enumerate}
\item[(i)] The vanishing
$\widetilde{\mathsf{K}}_{\alpha,\widehat{\alpha}}=0$ holds {unless}
$|{\alpha}|\geq |\widehat{\alpha}|$.
\item[(ii)]
The coefficients of  $\widetilde{\mathsf{K}}_{\alpha,\widehat{\alpha}}\in
\mathbb{Q}[i,c_1,c_2,c_3]((u))$ as a series in $u$ 
are homogeneous{\footnote{The variable $c_k$ has degree $k$
for the homogeneity (and the complex number $i$ has degree $0$).}} in the variables $c_i$
of degree $$|\alpha|+\ell(\alpha) - |\widehat{\alpha}| 
- \ell(\widehat{\alpha})-3(\ell(\alpha)-1).$$ 
\end{enumerate}
Let $\mathcal{T}_\nu\rightarrow \mathcal{X}$ be the rank 3 relative tangent bundle of the family
$\nu:\mathcal{X}\rightarrow \mathcal{Y}$.
Via the substitution
\begin{equation*} 
c_i=c_i(\mathcal{T}_\nu),
\end{equation*}
the matrix elements  of $\widetilde{\mathsf{K}}$
act by cup product on the cohomology 
of $\mathcal{X}$ with $\mathbb{Q}[i]((u))$-coefficients.

The matrix $\widetilde{\mathsf{K}}$ is 
used to define a correspondence
rule
\begin{equation}\label{pddff}
{\tau_{\alpha_1-1}(\gamma_1)\cdots
\tau_{\alpha_{\ell}-1}(\gamma_{\ell})}\ \  \mapsto\ \ 
\overline{\tau_{\alpha_1-1}(\gamma_1)\cdots
\tau_{\alpha_{\ell}-1}(\gamma_{\ell})}\ 
\end{equation}
with $\gamma_i\in H^*(\mathcal{X})$.
The left side is a product of standard descendent insertions.

The definition of the right side
of \eqref{pddff} requires a sum over all set
partitions $P$ of $\{ 1,\ldots, \ell \}$.
 For such a  set partition
$P$, each element $S\in P$
is a subset of $\{1,\ldots, \ell\}$.
Let $\alpha_S$ be the associated subpartition of
$\alpha$, and let
$$\gamma_S = \prod_{j\in S}\gamma_j.$$
In case all cohomology classes $\gamma_i$ are even{\footnote{If classes are
    odd, a sign must be introduced in the correspondence rule.
Signs will be discussed in Section \ref{sgns}.  }}
we define the right side of the correspondence rule  \eqref{pddff} 
by
\begin{equation}\label{mqq23}
\overline{\tau_{\alpha_1-1}(\gamma_1)\cdots
\tau_{\alpha_{\ell}-1}(\gamma_{\ell})}
=
\sum_{P \text{ set partition of }\{1,\ldots,\ell\}}\ \prod_{S\in P}\ \sum_{\widehat{\alpha}}\tau_{\widehat{\alpha}}(\widetilde{\mathsf{K}}_{\alpha_S,\widehat{\alpha}}\cdot\gamma_S) \ .
\end{equation}
The second sum  in \eqref{mqq23} is over 
all partitions $\widehat{\alpha}$ of positive size. However, by the vanishing
of property (i),
$$\widetilde{\mathsf{K}}_{\alpha_S,\widehat{\alpha}}=0 \ \ \text{unless}
\ \ |{\alpha_S}|\geq |\widehat{\alpha}|\, , $$
the summation index  may be restricted to partitions $\widehat{\alpha}$ of positive size bounded
by $|\alpha_S|$.

The terms on the right side of \eqref{mqq23} are diagonal descendent
insertions \eqref{ff33}. For example, let $\ell=5$
and consider the 
set partition $P$ of $\{1,2,3,4,5\}$ defined by the data
$$S_1\cup S_2=\{1,2,3,4,5\}\, , \ S_1=\{1,2,3\}\, , \ S_2=\{4,5\}$$
Terms on the right side of \eqref{mqq23}
associated to $P$
are of the form
$$\tau_{\widehat{\alpha}^1}(\widetilde{\mathsf{K}}_{\alpha_{S_1},\widehat{\alpha}^1}\cdot \gamma_{S_1})\cdot \tau_{\widehat{\alpha}^2}(\widetilde{\mathsf{K}}_{\alpha_{{S_2}},\widehat{\alpha}^2}\cdot \gamma_{S_2})\, .
$$
Here, $\widehat{\alpha}^1$ and $\widehat{\alpha}^2$ are
two partitions of positive size -- the superscript index serves
to distinguish them as both are summed over in \eqref{mqq23}.


Suppose  $|\alpha_S|=|\widehat{\alpha}|$
in the second sum in \eqref{mqq23}.
The homogeneity property (ii)
then places a strong constraint. The $u$ coefficients of  
\begin{equation*}
\widetilde{\mathsf{K}}_{\alpha_S,\widehat{\alpha}}\in
\mathbb{Q}[i,c_1,c_2,c_3]((u))
\end{equation*}
are homogeneous
of degree 
\begin{equation}\label{d2399}
3-2\ell(\alpha_S)  
- \ell(\widehat{\alpha}) \, .
\end{equation}
For the matrix element $\widetilde{\mathsf{K}}_{\alpha_S,\widehat{\alpha}}$ to be nonzero, the degree 
\eqref{d2399}
must be non-negative. Since the lengths of $\alpha_S$ and 
$\widehat{\alpha}$ are at least 1,
non-negativity of \eqref{d2399} is only possible if 
$$\ell(\alpha_S)= \ell(\widehat{\alpha})=1\, .$$
Then, we also have $\alpha_S=\widehat{\alpha}$ since the sizes match.

The above argument shows that the descendents on the right side of  \eqref{mqq23}
all correspond to partitions of size {\em less} than $|\alpha|$
except for the {\em leading term} obtained from the 
the maximal set partition
$$\{1\} \cup \{2\} \cup \ldots \cup \{\ell\} = \{1,2,\ldots, \ell\}$$
in $\ell$ parts.
The leading term of the descendent correspondence, calculated
in \cite{PPDC}, is a third basic property of $\widetilde{\mathsf{K}}$:

\begin{enumerate}
\item[(iii)]
$\ \ \overline{\tau_{\alpha_1-1}(\gamma_1)\cdots
\tau_{\alpha_{\ell}-1}(\gamma_{\ell})}
= (iu)^{\ell(\alpha)-|\alpha|}\, \tau_{\alpha_1-1}(\gamma_1)\cdots
\tau_{\alpha_{\ell}-1}(\gamma_{\ell}) +\ldots .$
\end{enumerate}


In case $\alpha=1^\ell$ has all parts equal to 1, then $\alpha_S$ also
has all parts equal to 1 for every $S\in P$. By property (ii), the $u$ coefficients of  
$\widetilde{\mathsf{K}}_{\alpha_S,\widehat{\alpha}}$
are homogeneous
of degree 
$$3-\ell(\alpha_S) - |\widehat{\alpha}| 
- \ell(\widehat{\alpha}),$$
and hence vanish unless
$$\alpha_S= \widehat{\alpha}=(1)\ .$$
Therefore, 
if
$\alpha$ has all parts equal to $1$,
the leading term is therefore the entire formula. 
We obtain a fourth property of the matrix $\widetilde{\mathsf{K}}$:

\begin{enumerate}
\item[(iv)]
$\ \overline{\tau_{0}(\gamma_1)\cdots
\tau_{0}(\gamma_{\ell})}
=  \tau_{0}(\gamma_1)\cdots
\tau_{0}(\gamma_{\ell})\, .$
\end{enumerate}

The geometric construction of $\widetilde{\mathsf{K}}$ in \cite{PPDC} 
expresses the coefficients explicitly in terms of the 1-legged capped
descendent vertex for stable pairs and stable maps.  These vertices
can be computed (as a rational function in the stable pairs case and
term by term in the genus parameter for stable maps). Hence, the coefficient
$$\widetilde{\mathsf{K}}_{\alpha,\widehat{\alpha}}\in \mathbb{Q}[i,c_1,c_2,c_3]((u))$$
can, in principle, be calculated term by term in $u$. The calculations
in practice are quite difficult, and 
complete closed formulas are not known{\footnote{See \cite{MorOOP,oop}
    for formulas is the stationary case for fixed
    $3$-folds $X$.}} for all of the coefficients.

\subsection{Conjecture I}\label{descor}
To state the Gromov-Witten/Pairs correspondence proposed and studied in
\cite{PPDC,PPQ} for the family $\nu:\mathcal{X}\rightarrow \mathcal{Y}$, the basic degree
$$d_\beta = \int_{\beta} c_1(T_\nu) \ \in \mathbb{Z}$$
associated to the fiber class $\beta$ is required.

\vspace{8pt}
\noindent{\bf Conjecture I [GW/P descendent correspondence].}
{\em For $\gamma_i \in H^{*}(\mathcal{X})$, we have
\begin{multline*}
(-q)^{-d_\beta/2}\ZZ_{\mathsf{P}}\Big(\nu;q\ \Big|  
{\tau_{\alpha_1-1}(\gamma_1)\cdots
\tau_{\alpha_{\ell}-1}(\gamma_{\ell})}
\Big)_\beta \\ =
(-iu)^{d_\beta}\ZZ'_{\mathsf{GW}}\Big(\nu;u\ \Big|   
\ \overline{\tau_{\alpha_1-1}(\gamma_1)\cdots
\tau_{\alpha_{\ell}-1}(\gamma_{\ell})}\ 
\Big)_\beta 
\end{multline*}
under the variable change $-q=e^{iu}$.}

\vspace{10pt}

Since the stable pairs side of the correspondence
 $$\ZZ_{\mathsf{P}}\Big(\nu;q\ \Big|  
{\tau_{\alpha_1-1}(\gamma_1)\cdots
\tau_{\alpha_{\ell}-1}(\gamma_{\ell})}
\Big)_\beta\, \in H^*(\mathcal{Y})((q))$$
 is  defined as a series in $q$, the change of variable
$-q=e^{iu}$ is {\em not} a priori well-defined.
However,
the stable pairs descendent series
is conjectured to be a rational function in $q$.
The change of variable $-q=e^{iu}$ is well-defined for a rational function
in $q$ 
by substitution. The well-posedness of 
the descendent correspondence depends upon the rationality conjecture.

\section{Results: past and future}

\subsection{Toric $3$-folds}
The basic case of the $\sf{GW}/\sf{P}$ descendent correspondence 
for families
is for the torus equivariant
Gromov-Witten and stable pairs theories of a nonsingular projective
toric 3-fold $\mathsf{X}$.
The torus equivariant theory is the families
theory for the homotopy quotient
$$ \nu:\mathcal{X} \, =\,  \mathsf{X} \times_{\mathsf{T}} \mathsf{ET} \rightarrow
\mathsf{BT} \, =\,  \mathcal{Y}\, ,$$
where $\mathsf{T}$ is the torus and $\mathsf{BT}$ has
an algebraic approximation by products of projective spaces.
In case of primary insertions (where all descendents
are of the form $\tau_0(\gamma)$), the correspondence was
conjectured in \cite{MNOP1,MNOP2} and proven in \cite{moop}. The
descendent
correspondence in the toric case
was suggested in \cite{MNOP2,pt} and proven in \cite{PPDC}. See \cite{part1,PPstat,
  PP2,pt2,max}
for stable pairs calculations in the toric case.

\vspace{8pt}
\noindent{\bf Theorem 1 [P.-Pixton].} {\em The {\sf{GW/P}} descendent
  correspondence  holds for the torus equivariant theories of
  nonsingular projective toric $3$-folds $\mathsf{X}$.}
\vspace{8pt}

Since Virasoro constraints are proven for the
descendent Gromov-Witten theory of toric $3$-folds \cite{Giv,Tel}, Theorem 1
implies Virasoro constraints for moduli spaces of sheaves, see
\cite{BLM, Mor, MorOOP,VB} for developments in these directions.

\subsection{Relative theories}
In order to use Theorem 1 to prove the descendent
correspondence for more general
$3$-folds $X$, a degeneration
strategy was pursued in \cite{PPQ}. There are several steps, but the
end result is a proof for $3$-folds $X$ that admit certain simple
degenerations to toric $3$-folds. A basic case{\footnote{For a fixed
    Calabi-Yau 3-fold $X$, the descendent and primary  correspondence
    is equivalent since the moduli spaces are all essentially 0 dimensional.}}
was the
Calabi-Yau quintic $3$-fold $X_5 \subset \mathbb{CP}^4$.

\vspace{8pt}
\noindent{\bf Theorem 2 [P.-Pixton]}. {\em The {\sf{GW/P}} correspondence
  holds for the $3$-fold $X_5$}.
\vspace{8pt}

Crucial to the strategy of \cite{PPQ} is the lifting of 
the entire descendent correspondence to the situation of
relative Gromov-Witten and stable pairs theories for
families $\nu$ over a 1-dimensional base $\mathcal{Y}$
where the fibers are allowed to have normal crossing
degenerations into two components.
As a consequence, many interesting relative cases of the
correspondence were proven in \cite{PPQ}.
For example, let $S_4\subset \mathbb{CP}^3$
be a nonsingular $K3$ quartic surface.

\vspace{8pt}
\noindent{\bf Theorem 3 [P.-Pixton]}. {\em The {\sf{GW/P}} descendent
  correspondence
  holds for the log Calabi-Yau $3$-fold $(\mathbb{CP}^3/S_4)$}.
\vspace{8pt}

In the past few years, Maulik and Ranganathan have promoted
the descendent correspondence conjecture to the logarithmic case
where the family $\nu$ is allowed to have arbitrary normal crossings
degenerations \cite{MR}. Their method should yield further
cases of the GW/P descendent correspondence for more general 
$3$-folds $X$ and log $3$-folds $(X/Y)$.

\subsection{Non-negative geometries} In recent work \cite{Par}, Pardon
has used new transversality arguments to reduce the case of the
correspondence for primary insertions to the correspondence for
local curves proven in \cite{BP,OP}. The outcome is the following result.

\vspace{8pt}
\noindent{\bf Theorem 4 [Pardon].} {\em The {\sf{GW/P}}
correspondence holds for primary insertions for
families $\nu:\mathcal{X} \rightarrow \mathcal{Y}$
where the fibers are nonsingular projective Fano or Calabi-Yau 3-folds.}
\vspace{8pt}

In particular, Pardon's results yield a new proof of Theorem 2 and
a new approach to the {\sf GW/P} correspondence in the primary case  for
rich geometries such as $S_4 \times E$
where $E$ is an elliptic curve \cite{ObP}. The GW/P descendent
correspondence{\footnote{The $3$-folds here
    are noncompact, so residue theories are
    required for the definitions. As usual, care must be taken to
    treat the reduced virtual class properly.}}
for the universal family
$$\nu: \mathcal{X}\, =\,  \mathbb{C} \times
\mathcal{S} \rightarrow
\mathcal{M}^{K3}\, =\, \mathcal{Y}$$
over the moduli of $K3$ surfaces is an interesting open direction.

\subsection{Moduli of curves}

A beautiful family of varieties is the universal 
curve, $\mathcal{C} \rightarrow \overline{\mathcal{M}}_{g,n}$,
over the moduli space of stable curves. An associated
family of $3$-folds is
$$\nu:\mathcal{X} \, = \, \mathbb{C}^2\times \mathcal{C} \rightarrow
\overline{\mathcal{M}}_{g,n} \, = \, \mathcal{Y}\, .$$
The definition of Gromov-Witten and stable pairs theories
for $\nu$ requires residues and relative geometry, see \cite{PHHT}.

\vspace{8pt}
\noindent{\bf Theorem 5 [P.- H.-H. Tseng].} {\em The {\sf{GW/P}}
  correspondence (with no insertions)
  holds for all relative boundary conditions for
the family $\mathbb{C}^2\times \mathcal{C} \rightarrow
\overline{\mathcal{M}}_{g,n}$.}
\vspace{8pt}

The {\sf{GW/P}} correspondence of Theorem 5 implies (and is
equivalent to) the Crepant Resolution Conjecture 
relating
the Gromov-Witten theory in all genera{\footnote{For the
    genus 0 theory, see \cite{BG,OPH} and the local curve
    results of \cite{BP,OP}.}}of the Hilbert scheme of
points of $\mathbb{C}^2$ to the orbifold Gromov-Witten
theory of the symmetric product of $\mathbb{C}^2$.

\vspace{-15pt}

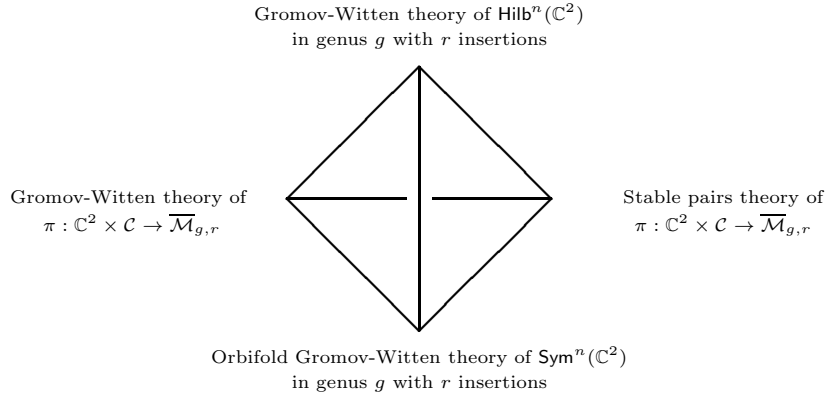
\begin{center}
\scriptsize
\begin{picture}(200,175)(-30,-50)
\thicklines
\put(25,25){\line(1,1){50}}
\put(25,25){\line(1,-1){50}}
\put(125,25){\line(-1,1){50}}
\put(125,25){\line(-1,-1){50}}
\put(75,-25){\line(0,1){100}}
\put(25,25){\line(1,0){45}}
\put(80,25){\line(1,0){45}}
\put(75,95){\makebox(0,0){Gromov-Witten theory of $\Hilb$}}
\put(75,85){\makebox(0,0){in genus $g$
with $r$ insertions}}
\put(75,-35){\makebox(0,0){Orbifold
Gromov-Witten theory of $\Sym(\com^2)$}}
\put(75,-45){\makebox(0,0){in genus $g$ 
with $r$ insertions}}
\put(190,25){\makebox(0,0){Stable pairs theory of}}
\put(190,15){\makebox(0,0){
$\pi:\mathbb{C}^2 \times 
\mathcal{C} \rightarrow
\overline{\mathcal{M}}_{g,r}$
}}
\put(-35,25){\makebox(0,0){Gromov-Witten theory of}}
\put(-35,15){\makebox(0,0){
$\pi:\mathbb{C}^2 \times 
\mathcal{C} \rightarrow
\overline{\mathcal{M}}_{g,r}
$}}
\end{picture}
\end{center}

\normalsize

\section{Stronger forms}
\subsection{Asymmetry for families}
The two sides of the {\sf{GW/P}} descendent correspondence
of Section \ref{descor}
are {\em not} symmetric for families
$$\nu:\mathcal{X} \rightarrow \mathcal{Y}\, .$$
The  
stable pairs side has a fully factored form
in terms of standard descendents
while the stable
maps side has diagonal descendents insertions \eqref{ff33}.
In case the base $\mathcal{Y}$ is a point, the K\"unneth
decomposition of the diagonal can be used to
express diagonal classes as sums of fully factorized terms,
but we do not have such K\"unneth decompositions over
arbitrary bases $\mathcal{Y}$.

\subsection{Conjecture II}
There is a canonical promotion of the {\sf{GW/P}} descendent correspondence
of Section \ref{descor} to include diagonal descendent insertions
\eqref{ff33} on both sides.

Let $\alpha=(\alpha_1, \ldots, \alpha_\ell)$ be a partition of positive size $|\alpha|$
and length $\ell$.
Let $D$ be a set partition of $\{1,\ldots, \ell\}$ in {\em ordered form}:
$$D_1 \cup \cdots \cup D_d =\{1,\ldots,\ell\}\, $$
and the elements of $D_i$ come before the elements of $D_j$ if $i<j$.
For example,
$$D_1\cup D_2 = \{1,2,3\}\, , \ D_1=\{1,2\}\, ,\  D_2=\{3\}$$
is in ordered form.

We will define a correspondence rule
\begin{equation*}
{\tau_{\alpha_{D_1}}(\gamma_1)\cdots
\tau_{\alpha_{D_d}}(\gamma_d)}\ \  \mapsto\ \ 
\overline{ \tau_{\alpha_{D_1}}(\gamma_1)\cdots
\tau_{\alpha_{D_d}}(\gamma_d) }\ ,
\end{equation*}
where $\alpha_{D_i}$ is the partition obtained
from $\alpha$ by selecting the parts with indices
$k_1,\ldots, k_{\ell(D_i)}$ in $D_i$,
and
$$\tau_{\alpha_{D_i}}(\gamma_i) = \tau_{\alpha_{k_1}-1,\ldots,
  \alpha_{k_{\ell(D_i)}}-1}(\gamma_i\cdot \Delta_{\ell(D_i)})\, $$
is the diagonal descendent insertion
with $\gamma_i\in H^*(\mathcal{X})$ following the notation of \eqref{ff33}.

Set partitions of $\{1,\ldots, \ell\}$ are partially ordered
by refinement.
Given another set partition $P$ of  $\{1,\ldots,\ell\}$,
let $D\wedge P$ denote the finest set partition for which $D$ and $P$
are both refinements.
A part $I\in D\wedge P$ is simultaneously{\footnote{For the parts of $D$ in
    $I$, we always view them in their natural order. Later, the ordering
  will be relevant for the sign.}}
a union of parts of
$D$ and a union of parts of $P$ (and no proper subset of $I$
has this property),
$$I = D_{i_1} \cup \cdots \cup D_{i_m}\, , \ \ \
I = S_{j_1} \cup \cdots \cup S_{j_n}\,.$$
Let $I\subset \{1,\ldots,\ell\}$ have $|I|$ elements.

In case all cohomology classes $\gamma_i$ are even{\footnote{If classes are
    odd, a sign must be introduced in the correspondence rule.
    Signs will be discussed in Section \ref{sgns}.}},
the right side of the correspondence rule is:
\begin{equation} \label{meme4}
\overline{ \tau_{\alpha_{D_1}}(\gamma_1)\cdots
\tau_{\alpha_{D_d}}(\gamma_d) }
=
\sum_{P \text{ set partition of }\{1,\ldots,\ell\}}\ \prod_{I\in  D\wedge P}\ \mathsf{T}_I \,,
\end{equation}
where we have
$$\mathsf{T}_I =\sum_{\widehat{\alpha}^1,\ldots, \widehat{\alpha}^n}
\tau_{\widehat{\alpha}^1\cup \cdots \cup \widehat{\alpha}^n}\left(
\prod_{k=1}^n\widetilde{\mathsf{K}}_{\alpha_{S_{j_k}},\widehat{\alpha}^k}\cdot
  \prod_{l=1}^m \gamma_{i_l}\cdot c_3(\mathcal{T}_\nu)^{|I|+1-n-m}\right)\, .
$$
Here, $\widehat{\alpha}^1\cup \cdots \cup \widehat{\alpha}^n$ is
the concatenation
of the partitions $\widehat{\alpha}^1,\ldots, \widehat{\alpha}^n$.

\vspace{8pt}
\noindent{\bf Conjecture II [GW/P diagonal descendent correspondence].}\\

\vspace{-4pt}
\noindent{\em For $\gamma_i \in H^{*}(\mathcal{X})$, we have
\begin{multline*}
(-q)^{-d_\beta/2}\ZZ_{\mathsf{P}}\Big(\nu;q\ \Big|  
{\tau_{\alpha_{D_1}}(\gamma_1)\cdots
\tau_{\alpha_{D_d}}(\gamma_d)}
\Big)_\beta \\ =
(-iu)^{d_\beta}\ZZ'_{\mathsf{GW}}\Big(\nu;u\ \Big|   
\
\overline{ \tau_{\alpha_{D_1}}(\gamma_1)\cdots
\tau_{\alpha_{D_d}}(\gamma_d) }
\ 
\Big)_\beta 
\end{multline*}
under the variable change $-q=e^{iu}$.}

\vspace{10pt}

As before, Conjecture II relies upon a rationality conjecture
for the diagonal descendent stable pairs theory of $\nu$.
Conjecture II specializes to Conjecture I when
$D$ is the set partition of $\{1,\ldots, \ell\}$ with every part
of size 1.
As a consequence of Conjecture II, the
diagonal descendent Gromov-Witten and stable pairs theories
of the family $\nu$ are equivalent.

\subsection{Signs} \label{sgns}
Signs play an important role both in Conjectures I and II.
When  the cohomology classes $\gamma_i\in H^*(\mathcal{X})$ are not all even,
a natural sign
must be included in correspondence rule \eqref{meme4}.

Given a 
set partition $P$ of $\{1,\ldots, \ell\}$ indexing 
the sum on the right side of the rule \eqref{meme4}, let 
$$I_1\cup \cdots \cup I_{\ell(D\wedge P)} = \{1,\ldots, \ell\}\, .$$
be the parts of $D\wedge P$.
The parts $I_j$ of $D\wedge P$ are unordered, but we choose an ordering.
We then 
obtain a permutation of the parts $D_1\cup \cdots \cup D_d$
by moving each part $D_i$ to the associated ordered part $D_i \subset I_j$  (and
respecting the original order of the $D_i$ in each $I_j$).
The permutation, in turn, determines a sign $\mathbf{s}(P)$
by the anti-commutation of the 
odd classes $\gamma_i$ associated to the parts $D_i$. We then write
\begin{equation*} 
\overline{ \tau_{\alpha_{D_1}}(\gamma_1)\cdots
\tau_{\alpha_{D_d}}(\gamma_d) }
=
\sum_{P \text{ set partition of }\{1,\ldots,\ell\}}\ (-1)^{\mathbf{s}(P)}
\prod_{j=1}^{\ell(D\wedge P)}\ \mathsf{T}_{I_j} \,.
\end{equation*}
The descendent 
$\overline{ \tau_{\alpha_{D_1}}(\gamma_1)\cdots
\tau_{\alpha_{D_d}}(\gamma_d) }$
is easily seen to have the
same commutation rules with respect to odd
cohomology as
$\tau_{\alpha_{D_1}}(\gamma_1)\cdots
\tau_{\alpha_{D_d}}(\gamma_d) $.

\subsection{Conjecture III}

What about the general descendent $\tau_{k_1,\ldots,k_r}(\delta)$
discussed in Sections \ref{prgw} and \ref{prdt}? 
Can we write a
{\sf{GW/P}} descendent correspondence for these descendents
in a symmetric form which generalizes Conjectures I and  II?
The answer is {\em yes}: there is a generalization
of Conjectures I and II for the general
descendent $\tau_{k_1,\ldots,k_r}(\delta)$, but the form
is more intricate.

Let $\alpha=(\alpha_1, \ldots, \alpha_\ell)$ be a partition of positive size $|\alpha|$
and length $\ell$.
We will define a correspondence rule
\begin{equation*}
{\tau_{\alpha_{1}-1, \ldots, \alpha_{\ell}-1}(\delta)}\ \  \mapsto\ \ 
\overline{ \tau_{\alpha_{1}-1, \ldots, \alpha_{\ell}-1}(\delta)}\ 
\end{equation*}
without any restrictions on the class $\delta\in H^*(\mathcal{X}^{\ell})$.

The right side of
the correspondence rule is
\begin{equation*} 
\overline{ \tau_{\alpha_{1}-1, \ldots, \alpha_{\ell}-1}(\delta)}\  =
\sum_{P \text{ set partition of }\{1,\ldots,\ell\}}\
\sum_{\widehat{\alpha}^1,\ldots, \widehat{\alpha}^{\ell(P)}}
\ \mathsf{U}_P( {\widehat{\alpha}^1,\ldots, \widehat{\alpha}^{\ell(P)}})  
\,.
\end{equation*}
The second sum is over all $\ell(P)$-tuples of partitions of
positive size (but only finitely many terms will be nonzero).

To define the contribution $\mathsf{U}_P({\widehat{\alpha}^1,\ldots, \widehat{\alpha}^{\ell(P)}})$, we first
order the parts of $P$,
$$S_1\cup \cdots\cup S_{\ell(P)}=\{1,\ldots, \ell\}$$
We place the elements within each $S_k$ in the
canonical increasing order (inherited from $\{1,\ldots, \ell\}$).
We thus obtain 
a canonical permutation $\theta_P\in \Sigma_\ell$.
The permutation $\theta_P$ yields a canonical automorphism
$$\theta_{P}: \mathcal{X}^\ell \rightarrow \mathcal{X}^\ell$$
over $\mathcal{Y}$ by permuting the product factors (to simplify the notation,
we use the same symbol for both the permutation and the automorphism).

For example, let $\ell=5$
and consider the 
set partition $P$ of $\{1,2,3,4,5\}$ defined by the data
$$S_1\cup S_2=\{1,2,3,4,5\}\, , \ S_1=\{2,4,5\}\, , \ S_2=\{1,3\}\, .$$
Then, the permutation $\theta_P\in \Sigma_5$ is
$$ 1\to 2\, ,\ \  2\to 4\, ,\ \  3\to 5\, , \ \ 4\to 1\, ,\ \ 5\to 3\, .$$

The contribution $\mathsf{U}_P({\widehat{\alpha}^1,\ldots, \widehat{\alpha}^{\ell(P)}})$ is a  general descendent insertion with
$\sum_{k=1}^{\ell(P)} \ell(\widehat{\alpha}^k)$
descendent indices,
$$\mathsf{U}_P({\widehat{\alpha}^1,\ldots, \widehat{\alpha}^{\ell(P)}})=
\tau_{\widehat{\alpha}^1_1-1,\ldots,\widehat{\alpha}^1_{\ell(\widehat{\alpha}^1)}-1,
  \ldots,
  \widehat{\alpha}^{\ell(P)}_1-1,\ldots,
  \widehat{\alpha}^{\ell(P)}_{\ell(\widehat{\alpha}^{\ell(P)})}-1}
(\delta_P({\widehat{\alpha}^1,\ldots, \widehat{\alpha}^{\ell(P)}}))\, .
$$
The last step is to define the class
$\delta_P({\widehat{\alpha}^1,\ldots, \widehat{\alpha}^{\ell(P)}})
\in H^*(\mathcal{X}^{\sum_{k=1}^{\ell(P)} \ell(\widehat{\alpha}^k)})$.

Let $\iota:\Delta_P\hookrightarrow \mathcal{X}^\ell$ be the product of diagonals,
$$\Delta_P = \prod_{k=1}^{\ell(P)} \text{pr}^*_{S_k}(\Delta_{S_k})\, ,$$
where $\Delta_{S_k}\subset \mathcal{X}^{S_k}$ is the small
diagonal of the factor
$\text{pr}_{S_k}: \mathcal{X}^\ell \to \mathcal{X}^{S_k}$
corresponding to $S_k\subset \{1,\ldots, \ell\}$.
There is a corresponding product of diagonals
$$ \widehat{\iota}:\Delta_{\widehat{\alpha}^1,\ldots, \widehat{\alpha}^{\ell(P)}}
\hookrightarrow \mathcal{X}^{\sum_{k=1}^{\ell(P)} \ell(\widehat{\alpha}^k)}$$
defined by pulling-back the small diagonals of
the product factors,
$$\Delta_{\widehat{\alpha}^1,\ldots, \widehat{\alpha}^{\ell(P)}}
= \prod_{k=1}^{\ell(P)} \text{pr}^*_{\widehat{\alpha}^k}(\Delta_{\widehat{\alpha}^k})\, ,$$
where
$\text{pr}_{\widehat{\alpha}^k}: \mathcal{X}^{\sum_{k=1}^{\ell(P)} \ell(\widehat{\alpha}^k)} \to \mathcal{X}^{\ell(\widehat{\alpha}^k)}$
corresponds to the indices of $\widehat{\alpha}^k$. There is
a canonical isomorphism
$$\phi_P: \Delta_P \to \Delta_{\widehat{\alpha}^1,\ldots, \widehat{\alpha}^{\ell(P)}}$$
which takes the factor $\mathcal{X}$ corresponding to $S_k$ on the domain to
the factor $\mathcal{X}$ corresponding to $\widehat{\alpha}^k$ on the target.
Then,
$$
\delta_P({\widehat{\alpha}^1,\ldots, \widehat{\alpha}^{\ell(P)}})
=
\widehat{\iota}_*\left(
\prod_{k=1}^{\ell(P)}\widetilde{\mathsf{K}}_{\alpha_{S_k},\widehat{\alpha}^k}\cdot
\phi_{P*}\big(\iota^*
{\theta}_{P}^*(\delta)
\big)\right)\, ,
$$
where each term $\widetilde{\mathsf{K}}_{\alpha_{S_k},\widehat{\alpha}^k}$
acts on the factor $\mathcal{X}$ corresponding to $\widehat{\alpha}^k$ of 
$\Delta_{\widehat{\alpha}^1,\ldots, \widehat{\alpha}^{\ell(P)}}$.

\vspace{8pt}
\noindent{\bf Conjecture III [GW/P families descendent correspondence].}\\

\vspace{-4pt}
\noindent{\em For $\delta \in H^*(\mathcal{X}^\ell)$, we have
\begin{multline*}
(-q)^{-d_\beta/2}\ZZ_{\mathsf{P}}\Big(\nu;q\ \Big|  
{\tau_{\alpha_1-1,\ldots,\alpha_\ell-1}(\delta)}
\Big)_\beta \\ =
(-iu)^{d_\beta}\ZZ'_{\mathsf{GW}}\Big(\nu;u\ \Big|   
\
\overline{\tau_{\alpha_1-1,\ldots,\alpha_\ell-1}(\delta)}
\ 
\Big)_\beta 
\end{multline*}
under the variable change $-q=e^{iu}$.}

\vspace{10pt}

Conjecture III relies upon a rationality conjecture
for the general descendent stable pairs theory of $\nu$.
A nice exercise is to specialize Conjecture III to
the cases of Conjectures I and II and to derive
the sign rules there from Conjecture III.

At the moment, Conjectures II and III are known for families
only in
cases where K\"unneth decompositions are available and
Conjecture I is proven. Perhaps the most interesting
known case is that of the $\mathsf{T}$-equivariant theory
of a nonsingular projective toric $3$-fold $\mathsf{X}$.

\section{Acknowledgements}
The formulation of the {\sf{GW/P}} descendent correspondence and several related
results are from
joint work \cite{PPDC,PPQ}  with A. Pixton. Discussions and
collaborations with J. Bryan, C. Faber, D. Maulik, M. Moreira, G. Oberdieck,
A. Oblomkov,
A. Okounkov, N. Nekrasov, D. Ranganathan, M. Schimpf,
and H.-H. Tseng have played an important role. I would like to
thank J. Pardon for conversations related to Theorem 4.
Thanks also to
the ICBS and the Beijing Institute of Mathematical Sciences and
Applications for the opportunity to speak about the
descendent correspondence in July 2024 (and to revisit the correspondence in
families).
I was supported by SNF-200020-219369 and SwissMAP.

\address{Department of Mathematics, ETH Z\"urich\\
 Z\"urich, Switzerland.\\
\email{rahul@math.ethz.ch}}


\begin{thebibliography}{MNOP2}

\bibitem{BehFan}
K.~Behrend and B.~Fantechi,
\newblock {\em The intrinsic normal cone,} 
{Invent. Math.} {\bf 128} (1997), 45--88.

\bibitem{BLM} A.~Bojko, W.~Lim, and M.~Moreira,
  {\em Virasoro constraints for moduli of sheaves and vertex algebras},
  Invent. Math. {\bf 236} (2024), 387--476.


\bibitem{BG} J.~Bryan and T.~Graber, {\em The crepant resolution conjecture},
  In: Algebraic Geometry–Seattle 2005, Part 1, 23--42, Proc.
Sympos. Pure Math. 80, Amer. Math. Soc., Providence, RI, 2009.

\bibitem{BP} J.~Bryan and R.~Pandharipande,
{\em Local Gromov-Witten theory of curves}, JAMS {\bf{21}} (2008), 101--136.

\bibitem{Ast} A.~Douady and J.-L.~Verdier, {\em S\'eminaire
de G\'eom\'etrie Analytique de l'Ecole Normale Sup\'erieure
1974/1975}, Ast\'erisque {\bf 36-37} (1976),
Expos\'es VI-IX.

\bibitem{Giv} A.~Givental, {\em Gromov-Witten invariants
    and quantization of quadratic Hamiltonians}, Moscow Math. J. {\bf 1}
  (2001), 487--518.


\bibitem{LePot} J.~Le Potier, {\em Faisceaux semi-stable
et syst\`emes coh\'erents}, in {\em Vector bundles in algebraic
geometry (Durham, 1993)}, LMS Lecture Note Ser. {\bf 208}, 179--239,
Cambridge Univ. Press: Cambridge, 1995.


\bibitem{LiTian}
J.~Li and G.~Tian,
\newblock {\em 
Virtual moduli cycles and {G}romov-{W}itten invariants of algebraic
  varieties,} {JAMS}  {\bf 11}, 119--174, 1998.



  
\bibitem{MNOP1}
D.~Maulik, N.~Nekrasov, A.~Okounkov, and R.~Pandharipande,
\newblock {\em Gromov-{W}itten theory and {D}onaldson-{T}homas theory. {I}},
  Compos. Math. {\bf 142} (2006), 1263--1285.


\bibitem{MNOP2}
D.~Maulik, N.~Nekrasov, A.~Okounkov, and R.~Pandharipande,
\newblock {\em Gromov-{W}itten theory and {D}onaldson-{T}homas theory. {II}},
  Compos. Math. {\bf 142} (2006), 1286--1304.


\bibitem{moop}
D.~Maulik, A. ~Oblomkov, A.~Okounkov, and R.~Pandharipande,
\newblock {\em The Gromov-{W}itten/{D}onaldson-{T}homas correspondence
for toric 3-folds}, Invent. Math. {\bf 186} (2011), 435--479.

\bibitem{MR}
  D.~Maulik and D.~Ranganathan, {\em Logarithmic enumerative
    geometry for curves and sheaves}, arXiv:2311.14150 (2023).

\bibitem{Mor} M.~Moreira, {\em Virasoro conjecture for the
    stable pairs descendent theory of simply
    connected 3-folds (with applications
    to the Hilbert scheme of points of a surface)},
    Journal of the LMS {\bf 106} (2020), 154--191.


 

    
  \bibitem{MorOOP}
  M.~Moreira, A.~Oblomkov, A.~Okounkov, and R.~Pandharipande, {\em Virasoro constraints
    for stable pairs on toric 3-folds}, Forum of Mathematics Pi {\bf{10}} (2022).

   \bibitem{ObP} G.~Oberdieck and A.~Pixton, {\em Holomorphic
      anomaly equations and the Igusa cusp form conjecture}, Invent.
    Math. {\bf 213} (2018), 507--587.

  
\bibitem{oop} A.~Oblomkov, A.~Okounkov, and R.~Pandharipande, {\em GW/PT
descendent correspondence via vertex operators},
Comm. Math. Phys. {\bf 374} (2020), 1321--1359.


\bibitem{OPH} A.~Okounkov and R.~Pandharipande,
{\em Quantum cohomology 
of the Hilbert scheme of points of the plane}, Invent. Math. {\bf 179} (2010),
523--557.



\bibitem{OP} A.~Okounkov and R.~Pandharipande,
 {\em The local Donaldson-Tho\-mas theory of curves}, Geom. Topol. {\bf 14} (2010), 1503--1567.

\bibitem{PanDesc} R. Pandharipande, {\em Descendents for stable pairs
on 3-folds}, Modern Geometry: A celebration of the work of 
Simon Donaldson, Proc. Sympos. Pure Math. {\bf 99} (2018), 251--288.

  
\bibitem{part1}
R.~Pandharipande and A.~Pixton,
\newblock {\em Descendents on local curves: Rationality},  Comp. Math. {\bf 149} (2013), 81--124. 


\bibitem{PPstat}
R.~Pandharipande and A.~Pixton,
\newblock {\em Descendents on local curves: Stationary theory} 
in {\em Geometry and arithmetic}, 283--307, EMS Ser. Congr. Rep., Eur. Math. Soc., Z\"urich, 2012.


\bibitem{PP2}
R.~Pandharipande and A.~Pixton,
\newblock {\em Descendent theory for stable pairs on
toric 3-folds},  Jour. Math. Soc. Japan. {\bf 65} (2013),
1337--1372.


\bibitem{PPDC}
R.~Pandharipande and A.~Pixton,
\newblock {\em Gromov-Witten/Pairs descendent correspondence
for toric 3-folds}, Geom. Topol. {\bf 18} (2014), 2747--2821.

\bibitem{PPQ}
R.~Pandharipande and A.~Pixton,
\newblock {\em Gromov-Witten/Pairs correspondence
for the quintic}, JAMS {\bf 30} (2017), 389--449.

\bibitem{pt}
R.~Pandharipande and R.~P. Thomas,
\newblock {\em Curve counting via stable pairs in the derived
category}, Invent Math. {\bf 178} (2009), 407--447.

\bibitem{pt2}
R.~Pandharipande and R.~P. Thomas,
\newblock {\em The 3-fold vertex via stable pairs}, Geom. Topol.
{\bf 13} (2009), 1835--1876.

\bibitem{13H} R.~Pandharipande and R.~Thomas, {\em 13/2 ways of counting curves}
in {\em  Moduli spaces}, 282--333, London Math. Soc. Lecture Note Ser., 411, Cambridge Univ. Press, Cambridge, 2014.



\bibitem{PHHT} R.~Pandharipande and H.-H.~Tseng,
{\em Higher genus Gromov-Witten theory of $\text{Hilb}^n({\mathbb{C}}^2)$ and
CohFTs associated to local curves}, Forum of Mathematics Pi {\bf 7} (2019).



  
\bibitem{Par} J.~Pardon, {\em Universally counting curves
in Calabi-Yau threefolds}, arXiv:2308.02948 (2023).

\bibitem{max} M.~Schimpf, {\em Stable pairs on local curves and
    Bethe roots}, in preparation.


\bibitem{Tel} C.~Teleman, {\em The structure of 2D semisimple field theories},
  Invent. Math. {\bf 188} (2012), 525--588.

\bibitem{VB} D.~van Bree, {\em Virasoro constraints for
    moduli spaces of sheaves on surfaces}, Forum of Mathematics Sigma {\bf 11}
  (2023), 1--35.
  
\end{thebibliography}
\end{document}